\documentclass[12pt, reqno]{amsart}
\usepackage{ifthen,amsfonts,amsmath,amssymb,epic,eepic}
\usepackage{epsfig,color}

\makeatother

\theoremstyle{definition}

\def\fnum{equation}
\newtheorem{Thm}[\fnum]{Theorem}
\newtheorem{Cor}[\fnum]{Corollary}

\newtheorem{Lem}[\fnum]{Lemma}

\newtheorem{Def}[\fnum]{Definition}

\newtheorem{Pro}[\fnum]{Proposition}

\newcommand{\nn}{{\bf{n}}}

\newcommand{\diam}{{\text {diam}}}
\newcommand{\dist}{{\text {dist}}}

\newcommand{\sztp}{\Sigma^{0,2\pi}}

\def\ZZ{{\bold Z}}
\def\RR{{\bold R}}
\def\SS{{\bold S}}

\def\CC{{\bold C }}
\newcommand{\dv}{{\text {div}}}
\newcommand{\e}{{\text {e}}}
\newcommand{\Area}{{\text {Area}}}
\newcommand{\Genus}{{\text {gen}}}

\newcommand{\cP}{{\mathcal{P}}}

\newcommand{\ca}{{\mathcal{A}}}

\newcommand{\cB}{{\mathcal{B}}}

\newcommand{\cT}{{\mathcal{T}}}

\newcommand{\K}{{\text{K}}}
\newcommand{\vb}{V \negmedspace B}
\newcommand{\kg}{k_{g}}

\newcommand{\eqr}[1]{(\ref{#1})}

\begin{document}

\title[Multi-valued graphs]
{The space of embedded minimal surfaces of fixed genus in a $3$-manifold II;
Multi-valued graphs in disks}

\author{Tobias H. Colding}%
\address{Courant Institute of Mathematical Sciences and MIT\\
251 Mercer Street\\
New York, NY 10012 and 77 Mass. Av., Cambridge, MA 02139}
\author{William P. Minicozzi II}%
\address{Department of Mathematics\\
Johns Hopkins University\\
3400 N. Charles St.\\
Baltimore, MD 21218}
\thanks{The first author was partially supported by NSF Grant DMS 9803253
and an Alfred P. Sloan Research Fellowship
and the second author by NSF Grant DMS 9803144
and an Alfred P. Sloan Research Fellowship.}

\email{colding@cims.nyu.edu and minicozz@math.jhu.edu}

\maketitle

\numberwithin{equation}{section}

\section{Introduction} \label{s:s0}

This paper is the second in a series where we attempt to give a
complete description of the space of all embedded minimal surfaces of
fixed genus in a fixed (but arbitrary)
closed $3$-manifold. The key for understanding
such surfaces is to understand the local structure in a
ball and in particular the structure of an embedded minimal
disk in a ball in $\RR^3$.
We show here that if the curvature of such a disk
becomes large at some point, then it contains an almost flat
multi-valued graph nearby that continues almost all the way to the
boundary.

\begin{figure}[htbp]
    \setlength{\captionindent}{20pt}
    \begin{minipage}[t]{0.5\textwidth}
    \centering\input{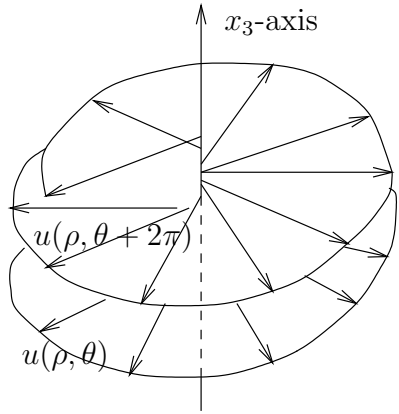}
    \caption{A multi-valued graph.} \label{f:1}
    \end{minipage}
\end{figure}

Let $\cP$ be the universal cover of the
punctured plane $\CC \setminus \{ 0 \}$ with global (polar)
coordinates $(\rho , \theta)$. An $N$-valued graph $\Sigma$ over
the annulus $D_{r_2} \setminus D_{r_1}$ (see fig. \ref{f:1}) is a (single-valued) graph
over
\begin{equation} \label{e:defmvg}
\{ (\rho ,\theta ) \in \cP \, | \, r_1 < \rho < r_2 {\text{ and }}
|\theta| \leq \pi \, N \} \, .
\end{equation}

\begin{Thm} \label{t:blowupwinding0}
Given $N\in \ZZ_+$, $\epsilon > 0$, there exist
$C_1,\,C_2>0$ so: Let
$0\in \Sigma^2\subset B_{R}\subset \RR^3$ be an embedded minimal
disk, $\partial \Sigma\subset \partial B_{R}$. If
$\max_{B_{r_0} \cap \Sigma}|A|^2\geq 4\,C_1^2\,r_0^{-2}$ for some
$R>r_0>0$, then there exists
(after a rotation)
an $N$-valued graph $\Sigma_g \subset \Sigma$ over $D_{R/C_2}
\setminus D_{2r_0}$ with gradient $\leq \epsilon$
and
$\Sigma_g \subset \{ x_3^2 \leq \epsilon^2 \, (x_1^2 + x_2^2) \}$.
\end{Thm}

This theorem is modeled by one half of the helicoid and its rescalings.
Recall that the helicoid is
the minimal surface $\Sigma^2$ in $\RR^3$
parameterized by
\begin{equation} \label{e:helicoid}
(s\,\cos t,s\sin t,t)
\end{equation}
where $s,t\in\RR$. By one half of the helicoid we mean the
multi-valued graph given by requiring that $s>0$ in \eqr{e:helicoid}.

Theorem \ref{t:blowupwinding0} will
follow by combining a blow up result with \cite{CM3}.
This blow up result says that if an embedded
minimal disk in a ball has large curvature at a point, then it
contains a small almost flat multi-valued graph nearby, that is:

\begin{Thm} \label{t:blowupwindinga}
See fig. \ref{f:2}.
Given $N , \omega>1$, and $\epsilon > 0$, there exists
$C=C(N,\omega,\epsilon)>0$ so: Let
$0\in \Sigma^2\subset B_{R}\subset \RR^3$ be an embedded minimal
disk, $\partial \Sigma\subset \partial B_{R}$. If
$\sup_{B_{r_0} \cap \Sigma}|A|^2\leq 4\,C^2\,r_0^{-2}$
and $|A|^2(0)=C^2\,r_0^{-2}$ for some $0<r_0<R$, then there exist
$ \bar{R} < r_0 / \omega$ and (after a rotation)
an $N$-valued graph $\Sigma_g \subset \Sigma$ over $D_{\omega \bar{R} }
\setminus D_{\bar{R} }$ with gradient $\leq \epsilon$, and
$\dist_{\Sigma}(0,\Sigma_g) \leq 4 \, \bar{R}$.
\end{Thm}

Recall that by the middle sheet $\Sigma^M$ of an $N$-valued graph $\Sigma$
we mean the portion over
\begin{equation}
\{ (\rho ,\theta ) \in \cP \, | \, r_1 < \rho < r_2 {\text{ and }}
0 \leq \theta \leq 2 \, \pi \} \, .
\end{equation}
The result that we need from \cite{CM3}
(combining  theorem 0.3 and lemma II.3.8 there) is:

\begin{figure}[htbp]
    \setlength{\captionindent}{20pt}
\begin{minipage}[t]{0.5\textwidth}
    \centering\input{shn3.pstex_t}
    \caption{Theorem \ref{t:blowupwindinga} -
finding a small multi-valued graph in
    a disk near a point of large curvature.} \label{f:2}
    \end{minipage}\begin{minipage}[t]{0.5\textwidth}
\centering\input{shn2.pstex_t}
    \caption{Theorem \ref{t:spin4ever2} - extending a small multi-valued graph
    in a disk.} \label{f:3}
\end{minipage}
\end{figure}

\begin{Thm} \label{t:spin4ever2}
\cite{CM3}; see fig. \ref{f:3}.
Given $N_1$ and $\tau > 0$, there exist $N , \Omega, \epsilon > 0$ so:
If $\Omega \, r_0 < 1 < R_0 / \Omega$, $\Sigma \subset B_{R_0}$ is
an embedded minimal disk with $\partial \Sigma \subset \partial B_{R_0}$,
and $\Sigma$ contains an $N$-valued minimal graph $\Sigma_g$ over
$D_1 \setminus D_{r_0}$ with gradient $\leq \epsilon$ and
$\Sigma_g \subset \{ x_3^2 \leq \epsilon^2 (x_1^2 + x_2^2) \}$, then
$\Sigma$ contains a $N_1$-valued graph $\Sigma_d$ over
$D_{R_0/\Omega} \setminus D_{r_0}$ with gradient $\leq \tau$ and
$(\Sigma_g)^M \subset \Sigma_d$.
\end{Thm}

As a consequence of Theorem \ref{t:blowupwinding0}, we will
show that if $|A|^2$ is blowing up for a sequence of embedded minimal disks,
then  there is a smooth minimal graph through this point in the limit of a
subsequence
(Theorem \ref{t:stablim} below).

Theorems \ref{t:blowupwinding0}, \ref{t:blowupwindinga},
\ref{t:spin4ever2}, \ref{t:stablim}
are local and are for simplicity
stated and proven only for $\RR^3$ with the flat metric although
they can with only very minor changes easily be seen to hold for
a sufficiently small ball in any
given fixed Riemannian $3$-manifold.

Let $x_1 , x_2 , x_3$ be the standard coordinates on $\RR^3$ and
$\Pi : \RR^3 \to \RR^2$ orthogonal projection to $\{ x_3 = 0 \}$.
For $y \in S \subset \Sigma \subset \RR^3$ and $s > 0$, the
extrinsic and intrinsic balls and tubes are
\begin{alignat}{2}
B_s(y) &= \{ x \in \RR^3 \, | \, |x-y| < s \} \, , \, & T_s(S) &=
\{ x \in \RR^3 \, | \, \dist_{\RR^3} (x , S) < s \} \, , \\
\cB_s(y) &= \{ x \in \Sigma \, | \, \dist_{\Sigma} (x , y) < s \}
\, , \, & \cT_s (S) &= \{ x \in \Sigma \, | \, \dist_{\Sigma} (x ,
S) < s \} \, .
\end{alignat}
$D_s$ denotes the disk $B_s(0) \cap \{ x_3 = 0 \}$.
$\K_{\Sigma}$ the sectional curvature of a smooth compact surface
$\Sigma$ and when
$\Sigma$ is immersed $A_{\Sigma}$ will be its second fundamental form.
When $\Sigma$ is oriented, $\nn_{\Sigma}$ is the unit normal.

\section{Poincar\'e and Caccioppoli
type inequalities for area and curvature}
\label{s:wtotcurv}

In this section, we will first estimate the area of a surface
(not necessarily minimal) in terms of its
total curvature; see Corollary \ref{c:poincare}.
This should be seen as
analogous to a Poincar\'e
inequality (for functions), and will be used similarly
later in this paper.
After that,
we will bound the curvature by the area for a minimal disk;
see Corollary \ref{c:lcacc}. This
inequality is
similar to a Caccioppoli inequality and, unlike the Poincar\'e type
inequality, relies on that the surface is minimal. Finally, we will apply
these inequalities to show a strengthened (intrinsic) version of a
result of Schoen and Simon.

\begin{Lem} \label{l:li1}
If $\cB_{r_0} (x) \subset \Sigma^2$
is disjoint from the cut locus of $x$,
\begin{align}
\text{Length}(\partial \cB_{r_0})-2\pi r_0&
=-\int_0^{r_0}\int_{\cB_{\rho}}\K_{\Sigma}\, ,
\label{e:oi1a}\\
\Area (\cB_{r_0}(x))- \pi \, r_0^2
&= - \int_{0}^{r_0} \int_0^{\tau} \int_{\cB_{\rho}(x)}
\K_{\Sigma}\, . \label{e:oi1}
\end{align}
\end{Lem}

\begin{proof}
For $0 < t \leq r_0$, by the Gauss-Bonnet
theorem,
\begin{equation} \label{e:oi2}
\frac{d}{dt} \int_{\partial \cB_t } 1
= \int_{\partial \cB_t} \kg =
2 \, \pi - \int_{\cB_t} \K_{\Sigma}
\, ,
\end{equation}
where $\kg$ is the geodesic curvature of $\partial \cB_t$.
Integrating \eqr{e:oi2} gives the lemma.
\end{proof}

\begin{Cor} \label{c:poincare}
If $\cB_{r_0} (x) \subset \Sigma^2$
is disjoint from the cut locus of $x$,
\begin{equation} \label{e:oi6}
\Area (\cB_{r_0}(x)) \leq \pi \, r_0^2 - \frac{1}{2} \, r_0^2 \,
\int_{\cB_{r_0}(x)} \min \{ \K_{\Sigma},0\} \, \, .
\end{equation}
\end{Cor}

\begin{Cor} \label{c:lcacc}
If $\Sigma^2\subset \RR^3$ is immersed and minimal,
$\cB_{r_0} \subset \Sigma^2$ is a disk, and
$\cB_{r_0} \cap \partial \Sigma = \emptyset$,
\begin{align}
t^2\int_{\cB_{{r_0}-2\,t}} |A|^2
&\leq r_0^{2}\int_{\cB_{r_0}}|A|^2\,(1-r/r_0)^2/2=
\int_{0}^{r_0} \int_0^{\tau} \int_{\cB_{\rho}(x)} |A|^2\notag\\
&=2\,(\Area \, (\cB_{r_0}) - \pi \, {r_0}^2) \leq r_0 \,
{\text{Length}} (\partial \cB_{r_0}) - 2 \pi  \, r_0^2
\, .\label{e:o1.23}
\end{align}
\end{Cor}

\begin{proof}
Since $\Sigma$ is minimal, $|A|^2=-2\,\K_{\Sigma}$ and
hence by Lemma \ref{l:li1}
\begin{equation} \label{e:bb1}
t^2 \, \int_{\cB_{{r_0}-2\,t}} |A|^2 \leq
t \, \int_{0}^{{r_0}-t} \int_{\cB_{\rho}} |A|^2 \leq
\int_{0}^{{r_0}} \int_0^{\tau} \int_{\cB_{\rho}} |A|^2 =
2 \, ( \Area \, (\cB_{r_0}) - \pi \, {r_0}^2 ) \, .
\end{equation}
The first equality follows by integration by parts twice
(using the coarea formula).
To get the last inequality in \eqr{e:o1.23}, note that
$\frac{d^2}{dt^2}{\text{Length}} (\partial \cB_{t})\geq 0$ by
\eqr{e:oi2} hence
$\frac{d}{dt}{\text{Length}} (\partial \cB_{t})\geq t\,
{\text{Length}} (\partial \cB_{t})$
and consequently $\frac{d}{dt}\left({\text{Length}}
(\partial \cB_{t}) /t\right)\geq 0$.
From this it follows easily.
\end{proof}

The following lemma and its corollary generalizes the main
result of \cite{ScSi}:

\begin{Lem} \label{l:scsi}
Given $C$, there exists $\epsilon>0$ so if
$\cB_{9 s}\subset \Sigma\subset \RR^3$ is an embedded minimal disk,
\begin{equation} \label{e:tcb0}
\int_{\cB_{9 s}}|A|^2 \leq C {\text{ and }}
\int_{\cB_{9 s}\setminus \cB_s} |A|^2 \leq \epsilon \, ,
\end{equation}
then $\sup_{\cB_{s}} |A|^2 \leq s^{-2}$.
\end{Lem}

\begin{proof}
Observe first that for $\epsilon$ small,  \cite{CiSc} and \eqr{e:tcb0} give
\begin{equation} \label{e:ciscann}
\sup_{\cB_{8 s}\setminus \cB_{2s}} |A|^2 \leq C_1^2 \,
\epsilon \, s^{-2} \, .
\end{equation}
By \eqr{e:oi1a} and \eqr{e:tcb0}
\begin{equation} \label{e:lbd}
{\text{Length}} (\partial \cB_{2s}) \leq
( 4 \pi + C) \, s \, .
\end{equation}
We will next use \eqr{e:ciscann} and \eqr{e:lbd} to show that,
after rotating $\RR^3$,  $\cB_{8 s}\setminus \cB_{2s}$ is
(locally) a graph over $\{ x_3 = 0 \}$ and furthermore $|\Pi
(\partial \cB_{8\,s})| > 3 \, s$.  Combining these two facts with
embeddedness, the lemma will then follow easily from Rado's
theorem.

By \eqr{e:lbd},  $\diam (\cB_{8s} \setminus \cB_{2s}) \leq (12  +
2 \pi + C/2)s$.  Hence, integrating \eqr{e:ciscann} gives
\begin{equation} \label{e:gmapaa}
\sup_{x,x' \in \cB_{8s} \setminus \cB_{2s}} \dist_{\SS^2} (\nn(x') , \nn(x) )
\leq C_1 \, \epsilon^{1/2} \, ( 12 + 2 \pi + C/2) \, .
\end{equation}
We can therefore rotate $\RR^3$ so that
\begin{equation} \label{e:gmapa}
\sup_{\cB_{8s} \setminus \cB_{2s}} |\nabla x_3|
\leq C_2 \, \epsilon^{1/2} \, ( 1 + C) \, .
\end{equation}

Given $y \in \partial \cB_{2s}$, let $\gamma_y$ be the outward normal
geodesic from $y$ to $\partial \cB_{8s}$ parametrized by arclength
on $[0,6s]$.    Integrating \eqr{e:ciscann} gives
\begin{equation} \label{e:ciscr1}
\int_{\gamma_y |_{[0,t]}} |k_g^{\RR^3}| \leq
\int_{\gamma_y |_{[0,t]} } |A| \leq C_1 \, \epsilon^{1/2} \, t / s \, ,
\end{equation}
where $k_g^{\RR^3}$ is the geodesic curvature of $\gamma_y$ in
$\RR^3$.
Combining \eqr{e:gmapa} with \eqr{e:ciscr1} gives (see fig. \ref{f:4})
\begin{equation} \label{e:gmapa3}
\langle \nabla |\Pi ( \cdot) - \Pi(y)| , \gamma_y' \rangle > 1 -
C_3 \, \epsilon^{1/2} \, ( 1 + C) \, .
\end{equation}
Integrating \eqr{e:gmapa3}, we get that for $\epsilon$ small,
$|\Pi (\partial \cB_{8\,s})| > 3 \, s$.

See fig. \ref{f:5}.
Combining $|\Pi (\partial \cB_{8\,s})| > 3 \, s$ and
\eqr{e:gmapa},  it follows that, for $\epsilon $ small,
$\Pi^{-1} (\partial D_{2s}) \cap
\cB_{8s}$ is a collection of immersed multi-valued graphs over
$\partial D_{2s}$.   Since $\cB_{8s}$ is embedded, $\Pi^{-1}
(\partial D_{2s}) \cap \cB_{8s}$ consists of disjoint embedded
circles which are graphs over $\partial D_{2s}$; this is the only
use of embeddedness.   Since $x_1^2+x_2^2$ is subharmonic on the
disk $\cB_{8\,s}$, these circles bound disks in $\cB_{8\,s}$ which
are then graphs by Rado's theorem (see, e.g., \cite{CM1}). The
lemma now follows easily from \eqr{e:gmapa} and the mean value
inequality.
\end{proof}

\begin{figure}[htbp]
    \setlength{\captionindent}{20pt}
\begin{minipage}[t]{0.5\textwidth}
\centering\input{blow3.pstex_t}
    \caption{Proof of Lemma \ref{l:scsi}:
By \eqr{e:gmapa} and
\eqr{e:ciscr1}, each $\gamma_y$ is almost a
horizontal line segment of length $6s$.  Therefore,
$|\Pi (\partial \cB_{8\,s})| > 3 \, s$.} \label{f:4}
\end{minipage}\begin{minipage}[t]{0.5\textwidth}
\centering\input{blow4.pstex_t}
    \caption{Proof of Lemma \ref{l:scsi}:
$\Pi^{-1} (\partial D_{2s})\cap \cB_{8s}$ is a union
of graphs over $\partial D_{2s}$.  Each bounds a
 graph in $\Sigma$ over
$D_{2s}$ by Rado's theorem.}
\label{f:5}
\end{minipage}
\end{figure}

\begin{Cor} \label{c:scsi}
Given $C_I$, there exists $C_P$ so
if $\cB_{2s}\subset \Sigma\subset \RR^3$ is an embedded
minimal disk with
\begin{equation} \label{e:tcbp}
\int_{\cB_{2s}} |A|^2 \leq C_I \, ,
\end{equation}
then $\sup_{\cB_{s}} |A|^2 \leq C_P \, s^{-2}$.
\end{Cor}

\begin{proof}
Let $\epsilon > 0$ be given by Lemma \ref{l:scsi} with $C=C_I$
and then let $N$ be the least integer greater than $C_I / \epsilon$.
Given $x \in \cB_s$, there exists $1\leq j \leq N$ with
\begin{equation}        \label{e:cth}
\int_{\cB_{9^{1-j} s}(x) \setminus \cB_{9^{-j} s}(x)} |A|^2 \leq
C_I / N \leq \epsilon \, .
\end{equation}
Combining \eqr{e:tcbp} and \eqr{e:cth}, Lemma \ref{l:scsi} gives that
$|A|^2(x) \leq (9^{-j} s)^{-2} \leq 9^{2N} \, s^{-2}$.
\end{proof}

We close this section with a generalization to
surfaces of higher genus; see Theorem \ref{t:scsig} below.
This will not be used in this paper
but will be useful in \cite{CM6}.  First we need:

\begin{figure}[htbp]
    \setlength{\captionindent}{20pt}
\begin{minipage}[t]{0.5\textwidth}
\centering\input{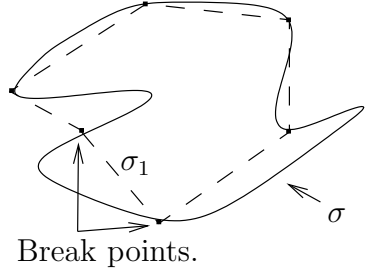}
    \caption{Lemma \ref{l:replacement}:
A curve $\sigma$ and a broken geodesic $\sigma_1$
homotopic to $\sigma$ in $\cT_{r_0}(\sigma)$.}
\label{f:6}
\end{minipage}
\end{figure}

\begin{Lem}   \label{l:replacement}
Let $\Sigma$ be a surface and $\sigma \subset \Sigma$ a
simple closed curve with length $\leq C\,r_0$.  If for all
$y\in \sigma$ the ball $\cB_{r_0}(y)$ is a disk disjoint
from $\partial \Sigma$, then
there is a broken geodesic
$\sigma_1 \subset \Sigma $ homotopic to
$\sigma$ in $\cT_{r_0}(\sigma)$
and with   $\leq C+1$ breaks; see fig. \ref{f:6}. If  $\Sigma$ is
an annulus with $\K_{\Sigma} \leq 0$ and $\sigma$ separates $\partial \Sigma$,
then $\sigma_1$ contains a simple curve $\sigma_2$ homotopic to $\sigma$ with
$\leq C+2$ breaks.
\end{Lem}

\begin{proof}
Parametrize $\sigma$ by arclength so that
$\sigma(0)=\sigma(\text{Length}(\sigma))$.
Let
$0=t_0<\cdots<t_n= \text{Length}(\sigma)$ be a subdivision with
$t_{i+1}-t_i \leq r_0$ and $n \leq C+1$.  Since $\cB_{r_0}(y)$
is a disk for all $y\in \sigma$, it follows that we can replace
$\sigma$ with a broken geodesic $\sigma_1$ with breaks at
$\sigma(t_i)=\sigma(t_i)$ and which is homotopic to $\sigma$
in $\cT_{r_0}(\sigma)$.

Suppose also now that $\Sigma$ is
an annulus with $\K_{\Sigma} \leq 0$ and $\sigma$ is topologically nontrivial.
Let $[a,b]$ be a maximal interval so that $\sigma_1 |_{[a,b]}$ is simple.
We are done if $\sigma_1 |_{[a,b]}$ is homotopic to $\sigma$.
Otherwise, $\sigma_1 |_{[a,b]}$ bounds a disk in $\Sigma$ and the Gauss-Bonnet
theorem implies that $\sigma_1 |_{(a,b)}$ contains a break.  Hence,
 replacing $\sigma_1$ by $\sigma_1 \setminus \sigma_1 |_{(a,b)}$
gives a subcurve homotopic to $\sigma$ but does not
increase the number of breaks.  Repeating
this eventually gives $\sigma_2$.
\end{proof}

Given a surface $\Sigma$ with boundary $\partial \Sigma$, we will
define the {\it genus} of $\Sigma$ ($\Genus (\Sigma)$) to be the
genus of the closed surface $\hat{\Sigma}$ obtained by adding a
disk to each boundary circle.  For example, the disk and the
annulus are both genus zero; on the other hand, a closed surface
of genus $g$ with $k$ disks removed has genus $g$.

In contrast to Corollary \ref{c:scsi} (and the results preceding it),
the next result concerns surfaces intersected with extrinsic balls.
Below,  $\Sigma_{0,s}$ is the component of $B_{s} \cap
\Sigma$ with $0 \in \Sigma_{0,s}$.

\begin{Thm} \label{t:scsig}
Given $C_a , g$, there exist $C_c , C_r$ so: If $0
\in \Sigma\subset B_{r_0}$ is an embedded minimal surface
with $\partial \Sigma\subset \partial B_{r_0}$, $\Genus(\Sigma)
\leq g$, $\Area (\Sigma) \leq C_a \, r_0^2$, and for each $C_r
\,r_0 \leq s \leq r_0$,
$\Sigma \setminus \Sigma_{0,s}$ is topologically an annulus,
then $\Sigma$ is a disk and
$\sup_{\Sigma_{0, C_r \, r_0}} |A|^2 \leq C_c \, r_0^{-2}$.
\end{Thm}

\begin{proof}
By the coarea formula, we can find  $r_0/2 \leq r_1 \leq 3r_0 / 4$
with ${\text{Length}}(\partial B_{r_1} \cap \Sigma) \leq 4 \, C_a \, r_0$.
It is easy to see from the maximum principle that
$\cB_{r_0/4}(y)$ is a disk for each $y \in
\partial B_{r_1} \cap \Sigma$ (we will take $C_r < 1/4$).
Applying Lemma \ref{l:replacement} to $\partial
\Sigma_{0,r_1} \subset \Sigma \setminus \Sigma_{0,r_0/4}$,
we get a simple broken geodesic
$\sigma_2 \subset \cT_{r_0/4}(\partial \Sigma_{0,r_1})$
homotopic to $\partial
\Sigma_{0,r_1}$  and with $\leq 16 \, C_a+2$ breaks.
Consequently, the Gauss-Bonnet
theorem  gives
\begin{equation}
        \int_{\Sigma_{0,r_0/4}} |A|^2
= - 2 \int_{\Sigma_{0,r_0/4}} K_{\Sigma} \leq
        8 \, \pi \, g + 2 \, \int_{\sigma_2} |k_g|
\leq 8 \, \pi ( g + 4\,  C_a + 1) \, .
\end{equation}
 For $\epsilon > 0$, arguing as in Corollary \ref{c:scsi} gives $r_2$ with
$\int_{\Sigma_{0,5r_2} \setminus \Sigma_{0,r_2}} |A|^2 \leq \epsilon^2$
so,     by \cite{CiSc},
\begin{equation}        \label{e:pcw}
        \sup_{\Sigma_{0,4r_2} \setminus \Sigma_{0,2r_2}} |A|^2
\leq C \, \epsilon^2 \, r_2^{-2} \, .
\end{equation}
 Using the area
bound,  $\partial \Sigma_{0,3r_2}$  can be covered by $C C_a$ intrinsic
balls $\cB_{r_2/4}(x_i)$ with  $x_i \in \partial \Sigma_{0,3r_2}$
(by the maximum principle, each $\cB_{r_2/4}(x_i)$ is a disk).
Hence, since  $\partial \Sigma_{0,3r_2}$ is connected,
any two points in $\partial \Sigma_{0,3r_2}$
can be joined
by a curve in $\Sigma_{0,4r_2} \setminus \Sigma_{0,2r_2}$ of length
$\leq C \,  r_2$.
Integrating \eqr{e:pcw} twice then gives a plane $P \subset \RR^3$
with $\partial \Sigma_{0,3r_2} \subset T_{C \, \epsilon \, r_2}(P)$.
By the convex hull property,
 $0 \in \Sigma_{0,3r_2} \subset T_{C \, \epsilon \, r_2}(P)$.
Hence, since $\partial \Sigma_{0,3r_2}$ is
connected and embedded, $\partial \Sigma_{0,3r_2}$ is a
graph over the boundary
of a convex domain for $\epsilon$ small.
The standard existence theory and Rado's theorem give
a minimal graph $\Sigma_g$ with $\partial \Sigma_g =
\partial \Sigma_{0,3r_2}$.  By translating $\Sigma_g$ above
$\Sigma_{0,3r_2}$ and sliding it down to the first point of
contact, and then repeating
this from below, it follows easily from the strong maximum
principle that $\Sigma_g =  \Sigma_{0,3r_2}$, completing the proof.
\end{proof}

\section{Finding large nearly stable pieces} \label{s:s1a}

We will collect here some results on stability of
minimal surfaces which will be used later to conclude that
certain sectors
are nearly stable.
The basic point is that
two disjoint but nearby embedded minimal surfaces satisfying a
priori curvature estimates must be nearly stable (made precise
below). We start by recalling the definition of $\delta_s$-stability.
Let again $\Sigma\subset \RR^3$ be an embedded oriented minimal surface.

\begin{Def} \label{d:almst}
($\delta_s$-stability).
Given $\delta_s \geq 0$, set
\begin{equation}
L_{\delta_s} = \Delta + (1-\delta_s) |A|^2 \, ,
\end{equation}
so that $L_0$ is the usual Jacobi operator on $\Sigma$. A domain
$\Omega \subset \Sigma$ is $\delta_s$-stable if $\int \phi \,
L_{\delta_s} \phi \leq 0$ for any compactly supported Lipschitz
function $\phi$ (i.e., $\phi \in C_0^{0,1}(\Omega)$).
\end{Def}

It follows that $\Omega$ is $\delta_s$-stable if and only if, for
all $\phi \in C_0^{0,1}(\Omega)$, we have the $\delta_s$-stability
inequality:
\begin{equation} \label{e:stabineq}
(1-\delta_s) \int |A|^2 \phi^2 \leq  \int |\nabla \phi |^2 \, .
\end{equation}

Since the Jacobi equation is the linearization of the minimal
graph equation over $\Sigma$, standard calculations give:

\begin{Lem} \label{l:mge}
There exists $\delta_g > 0$ so that if $\Sigma$ is minimal and $u
$ is a positive solution of the minimal graph equation over
$\Sigma$ (i.e., $\{ x + u(x) \, \nn_{\Sigma} (x) \, |\, x \in
\Sigma \}$ is minimal) with $|\nabla u| + |u| \, |A| \leq
\delta_g$, then $w = \log u$ satisfies on $\Sigma$
\begin{equation} \label{e:wlog2}
\Delta w = - |\nabla w|^2 + \dv (a \nabla w) + \langle \nabla w ,
a \nabla w \rangle + \langle b , \nabla w \rangle + (c-1) |A|^2 \,
,
\end{equation}
for functions $a_{ij} , b_j , c$ on $\Sigma$ with $|a| , |c| \leq
3 \, |A| \, |u| + |\nabla u|$ and $|b| \leq 2 \, |A| \, |\nabla
u|$.
\end{Lem}

The following slight modification of
a standard argument (see, e.g., proposition 1.26 of \cite{CM1}) gives
a useful sufficient condition for $\delta_s$-stability of a domain:

\begin{Lem} \label{l:fiscm}
There exists $\delta > 0$ so: If $\Sigma$ is minimal
and $u > 0$ is a solution of the minimal graph equation over
$\Omega \subset \Sigma$ with $|\nabla u| + |u| \, |A| \leq
\delta$, then
$\Omega$ is $1/2$-stable.
\end{Lem}

\begin{proof}
Set $w=\log u$ and choose a cutoff function
$\phi \in C_0^{0,1}(\Omega)$.
Applying Stokes' theorem
to $\dv ( \phi^2 \, \nabla w - \phi^2 \, a \, \nabla w)$,
substituting \eqr{e:wlog2}, and using
$|a| , |c| \leq 3 \, \delta ,
|b| \leq 2 \, \delta \, |\nabla w|$ gives
\begin{align} \label{e:lapw2}
(1 - 3 \, \delta) \, \int \phi^2 \, |A|^2
&\leq - \int \phi^2 \, |\nabla w|^2 + \int \phi^2 \,
\langle \nabla w , b + a \, \nabla w \rangle
+ 2 \int \phi \langle \nabla \phi , \nabla w
- a \, \nabla w \rangle \notag \\
&\leq (5 \delta - 1 ) \, \int
\phi^2 \, |\nabla w|^2 + 2 (1 + 3 \delta)
\int |\phi \, \nabla w| \, |\nabla \phi|
\, .
\end{align}
The lemma now follows easily from
the absorbing inequality.
\end{proof}

We will use Lemma \ref{l:fiscm} to see that disjoint embedded minimal
surfaces that are close are nearly stable (Corollary \ref{c:delst} below).
Integrating $|\nabla \dist_{\SS^2}( \nn (x) , \nn )| \leq
|A|$ on geodesics gives
\begin{equation} \label{e:gmap}
\sup_{x' \in \cB_s(x)} \dist_{\SS^2} (\nn(x') , \nn(x) )
\leq s \, \sup_{\cB_s(x)} |A| \, .
\end{equation}
By \eqr{e:gmap},
we can choose $0 < \rho_2 < 1/4$ so: If
$\cB_{2s}(x) \subset \Sigma$, $s \, \sup_{\cB_{2s}(x)} |A| \leq 4 \,
\rho_2$, and $t \leq s$, then the component
$\Sigma_{x,t}$ of $B_{t}(x) \cap \Sigma$ with $x \in \Sigma_{x,t}$
is a graph over $T_x \Sigma$ with gradient $\leq t/s$ and
\begin{equation} \label{e:gmap1}
\inf_{x' \in \cB_{2s}(x)} |x'-x| / \dist_{\Sigma}(x,x') > 9/10 \, .
\end{equation}
One consequence is that if
$t \leq s$ and we translate $T_x \Sigma$ so that $x \in T_x \Sigma$, then
\begin{equation} \label{e:gmap2}
\sup_{x' \in \cB_t(x)} \, | x' - T_x \Sigma|
\leq t^2 / s \, .
\end{equation}

\begin{Lem} \label{l:dnl}
There exist $C_0 , \rho_0 > 0$ so:
If $\rho_1 \leq \min \{\rho_0,\rho_2\}$,
$\Sigma_1 , \Sigma_2 \subset \RR^3$ are oriented
minimal surfaces, $|A|^2\leq 4$ on each $\Sigma_i$,
$x \in \Sigma_1 \setminus \cT_{4\rho_2}(\partial \Sigma_1)$, $y
\in B_{\rho_1}(x) \cap \Sigma_2 \setminus \cT_{4\rho_2}(\partial
\Sigma_2)$, and $\cB_{2 \rho_1}(x) \cap \cB_{2 \rho_1}(y) =
\emptyset$, then $\cB_{\rho_2}(y)$ is the graph $\{ z + u(z) \,
\nn (z) \}$ over a domain containing $\cB_{\rho_2/2}(x)$ with $u
\ne 0$ and $|\nabla u| + 4 \, |u| \leq C_0 \, \rho_1$.
\end{Lem}

\begin{proof}
Since $\rho_1 \leq \rho_2$,
\eqr{e:gmap1} implies that
$\cB_{ 2\rho_2}(x) \cap \cB_{ 2\rho_2}(y) = \emptyset$.
If $t\leq 9 \rho_2 / 5$, then $|A|^2\leq 4$ implies that the
components $\Sigma_{x,t} , \Sigma_{y,t}$ of $B_{t}(x) \cap
\Sigma_1 , B_{t}(y) \cap \Sigma_2$, respectively, with $x\in
\Sigma_{x,t} , y\in \Sigma_{y,t}$, are graphs with gradient $\leq
t/(2\rho_2)$ over $T_x \Sigma_1 , T_y \Sigma_2$ and have
$\Sigma_{x,t} \subset \cB_{ 2\rho_2}(x) , \Sigma_{y,t} \subset
\cB_{ 2\rho_2}(y)$. The last conclusion implies that
$\Sigma_{x,t} \cap \Sigma_{y,t} = \emptyset$.
It now follows that $\Sigma_{x,t} , \Sigma_{y,t}$ are
graphs over the same plane.  Namely,
if we set
$\theta =  \dist_{\SS^2} (\nn(x) , \{ \nn(y) , - \nn(y) \} )$, then
\eqr{e:gmap2}, $|x-y| < \rho_1$,
and $\Sigma_{x,t} \cap \Sigma_{y,t} = \emptyset$ imply that
\begin{equation} \label{e:nnclose1}
\rho_1 - (t/2 -\rho_1) \, \sin \theta + t^2 / (2 \rho_2) >
- t^2 / (2\rho_2 )\, .
\end{equation}
Hence, $\sin \theta <  \rho_1 / (t/2 -\rho_1) +  t^2 /
[(t/2 -\rho_1)\rho_2]$. For $\rho_0 / \rho_2$ small, $\cB_{\rho_2}(y)$ is a
graph with bounded gradient over $T_x \Sigma_1$. The lemma now
follows easily using the Harnack inequality.
\end{proof}

\begin{figure}[htbp]
    \setlength{\captionindent}{20pt}
    \begin{minipage}[t]{0.5\textwidth}
\centering\input{blow7.pstex_t}
    \caption{Corollary \ref{c:delst}:  Two sufficiently close disjoint
minimal surfaces with bounded curvatures must be nearly stable.}
\label{f:7}
\end{minipage}\begin{minipage}[t]{0.5\textwidth}
\centering\input{blow8.pstex_t}
    \caption{The set $VB$ in \eqr{e:bplus}. Here
$x \in VB$ and $y \Sigma\setminus VB$.}
\label{f:8}
\end{minipage}
\end{figure}

Combining   Lemmas \ref{l:fiscm} and \ref{l:dnl} gives:

\begin{Cor} \label{c:delst}
See fig. \ref{f:7}.  Given $C_0 , \, \delta > 0$,
there exists
$\epsilon (C_0 , \delta) > 0$ so that
if $p_i \in \Sigma_i \subset \RR^3$ ($i=1,2$) are embedded minimal surfaces,
$\Sigma_1 \cap \Sigma_2 = \emptyset$,
$\cB_{2 R}(p_i) \cap \partial \Sigma_i = \emptyset$,
$ |p_1 - p_2| < \epsilon \, R$, and
\begin{equation} \label{e:curvhinc}
\sup_{\cB_{2 R}(p_i)} |A|^2 \leq
C_0 \, R^{-2} \, ,
\end{equation}
then $\cB_{R} (\tilde{p}_i) \subset \tilde{\Sigma}_i$ is $\delta$-stable where
$\tilde{p}_i$ is the point over $p_i$ in the universal cover $\tilde{\Sigma}_i$
of $\Sigma_i$.
\end{Cor}

The next result gives a decomposition
of an embedded minimal surface with bounded curvature into a portion with
bounded area and a union of disjoint $1/2$-stable domains.

\begin{Lem} \label{l:deltstb}
There exists $C_1$ so:
If $0 \in \Sigma \subset B_{2R} \subset \RR^3$ is an
embedded minimal surface with
$\partial \Sigma \subset \partial B_{2R}$, and $|A|^2 \leq 4$, then
there exist disjoint
$1/2$-stable subdomains $\Omega_j \subset \Sigma$ and
a function $\chi \leq 1$ which vanishes on $B_R \cap \Sigma \setminus
\cup_j \Omega_j$ so that
\begin{align} \label{e:abd}
\Area ( \{ x \in B_R \cap \Sigma \, | \, \chi (x) < 1 \} ) &
\leq C_1 \, R^3 \, , \\
\int_{\cB_R} |\nabla \chi|^2 & \leq C_1 \, R^3 \, . \label{e:chi}
\end{align}
\end{Lem}

\begin{proof}
We can assume that $R > \rho_2$ (otherwise $B_R \cap \Sigma$ is stable).
Let $\delta > 0$ be from Lemma \ref{l:fiscm}
and $C_0 , \rho_0$ be from Lemma \ref{l:dnl}. Set
$\rho_1 = \min \{ \rho_0 / C_0 , \delta / C_0 , \rho_2 / 4 \}$.

Given $x \in B_{2R - \rho_1} \cap \Sigma$,
let $\Sigma_{x}$ be the component of $B_{\rho_1}(x)\cap \Sigma$
with $x \in \Sigma_{x}$ and let $B_x^{+}$ be the component
of $B_{\rho_1}(x) \setminus \Sigma_{x}$
which $\nn(x)$ points into.
See fig. \ref{f:8}.  Set
\begin{equation} \label{e:bplus}
\vb = \{ x \in B_R \cap \Sigma \, | \,
B_x^{+} \cap \Sigma \setminus
\cB_{4 \, \rho_1 }(x) = \emptyset \}
\end{equation}
and let $\{ \Omega_j \}$ be the
components of $B_{R} \cap \Sigma \setminus \overline{\vb}$.
Choose a maximal disjoint collection $\{ \cB_{\rho_1}(y_i) \}_{1 \leq i \leq \nu}$
of balls centered in
$\vb$. Hence, the union of the balls
$\{ \cB_{2 \, \rho_1}(y_i) \}_{1 \leq i \leq \nu}$ covers $\vb$.
Further,
the ``half-balls'' $B_{\rho_1/2}(y_i) \cap
B_{y_i}^{+} $ are pairwise disjoint.  To see this, suppose that
$|y_i - y_j| < \rho_1$ but $y_j \notin \cB_{2 \rho_1}(y_i)$.
Then, by \eqr{e:gmap1},
$y_j \notin \cB_{8 \rho_1}(y_i)$ so
$\cB_{4 \rho_1}(y_j) \notin  B_{y_i}^+$ and
$\cB_{4 \rho_1}(y_i) \notin  B_{y_j}^+$; the triangle inequality
then implies that
$B_{\rho_1/2}(y_i) \cap
B_{y_i}^{+} \cap B_{\rho_1/2}(y_j) \cap
B_{y_j}^{+} = \emptyset$ as claimed.
By  \eqr{e:gmap}--\eqr{e:gmap2}, each
$B_{\rho_1/2}(y_i) \cap
B_{y_i}^{+}$
has volume approximately $\rho_1^{3}$
and is contained in $B_{2R}$ so that
$\nu \leq C \, R^3$.
Define the function $\chi$ on $\Sigma$ by
\begin{equation}
\chi (x)=
\begin{cases}
0 & \hbox{ if } x \in \vb \, , \\
\dist_{\Sigma} (x, \vb ) / \rho_1
& \hbox{ if } x\in \cT_{\rho_1}(\vb)\setminus \vb \, , \\
1 & \hbox{ otherwise }\, .
\end{cases}
\end{equation}
Since $\cT_{\rho_1}(\vb) \subset
\cup_{i=1}^{\nu}
\cB_{3 \, \rho_1}(y_i)$, $|A|^2 \leq 4$, and $\nu \leq C \, R^3$,
we get \eqr{e:abd}. Combining \eqr{e:abd} and $|\nabla \chi| \leq \rho_1^{-1}$
gives \eqr{e:chi} (taking $C_1$ larger).

It remains to show that each $\Omega_j$ is $1/2$-stable.
Fix $j$. By construction, if $x \in \Omega_j$, then
there exists $y_x \in B_x^{+} \cap \Sigma \setminus
\cB_{4 \, \rho_1 }(x)$ minimizing $|x - y_x|$ in
$B_x^{+} \cap \Sigma$.
In particular, by Lemma \ref{l:dnl},
$\cB_{\rho_2}(y_x)$ is the graph
$\{ z + u_{x}(z) \, \nn (z) \}$ over a domain containing $\cB_{\rho_2/2}(x)$
with $u_x > 0$ and
$|\nabla u_x| + 4 \, |u_x| \leq \min \{
\delta , \rho_0 \}$.
Choose a maximal disjoint collection of balls $\cB_{\rho_2 / 6}(x_i)$
with $x_i \in \Omega_j$ and let $u_{x_i} > 0$ be the corresponding
functions defined on $\cB_{\rho_2 / 2}(x_i)$.
Since $\Sigma$ is embedded (and compact) and $|u_{x_i}| < \rho_0$,
Lemma \ref{l:dnl} implies that
$u_{x_i}(x) = \min_{t>0} \{ x + t \, \nn(x) \in \Sigma \}$ for
 $x \in \cB_{\rho_2 / 2}(x_{i})$.
Hence,  $u_{x_{i_1}}(x) = u_{x_{i_2}}(x)$ for
$x \in \cB_{\rho_2 / 2}(x_{i_1}) \cap
\cB_{\rho_2 / 2}(x_{i_2})$.
Note that
$\cT_{\rho_2 / 6} (\Omega_j) \subset \cup_i \cB_{\rho_2 / 2}(x_i)$.
We conclude that the $u_{x_i}$'s give a well-defined function
$u_j> 0$ on $\cT_{\rho_2/6}(\Omega_j)$ with $|\nabla u_j| + |u_j| \, |A| \leq
\delta$.
Finally,
Lemma \ref{l:fiscm} implies that each $\Omega_j$ is
$1/2$-stable.
\end{proof}

\section{Total curvature and area of embedded minimal disks} \label{s:s1}

Using the decomposition of Lemma \ref{l:deltstb}, we next
obtain polynomial bounds
for the area and total curvature of intrinsic balls
in embedded minimal
disks with bounded curvature.

\begin{Lem} \label{l:vgr}
There exists $C_1$ so
if $0 \in \Sigma \subset B_{2\, R}$ is an
embedded minimal
disk, $\partial \Sigma \subset \partial B_{2\, R}$,
$|A|^2 \leq 4$, then
\begin{equation} \label{e:tcgr}
 \int_{0}^R
\int_{0}^t \int_{\cB_s} |A|^2 \, ds   \, dt =
2( \Area (\cB_{R})- \pi \, R^2)
\leq 6 \,\pi \, R^2 + 20 \, C_1 \, R^5 \,  .
\end{equation}
\end{Lem}

\begin{proof}
Let $C_1$, $\chi$, and $\cup_j \Omega_j$ be given by
Lemma \ref{l:deltstb}.
Define $\psi$ on $\cB_R$
by $\psi =\psi (\dist_{\Sigma}(0, \cdot) ) =
1- \dist_{\Sigma}(0, \cdot) /R$, so $\chi \psi$
vanishes off of $\cup_j \Omega_j$.
Using $\chi \psi$ in the $1/2$-stability inequality, the
absorbing inequality and \eqr{e:chi}
give
\begin{align} \label{e:stab1}
\int |A|^2 \chi^2 \psi^2 & \leq 2 \, \int |\nabla ( \chi \psi)|^2 =
2 \, \int \, \left( \chi^2 |\nabla \psi|^2 + 2 \chi \, \psi \langle
\nabla \chi , \nabla \psi \rangle + \psi^2 |\nabla \chi|^2 \right)
\notag \\ & \leq 6 \, C_1 \, R^3 + 3 \int \chi^2 |\nabla \psi|^2
\leq 6 \, C_1 \, R^3 + 3\,R^{-2} \Area (\cB_R)\, .
\end{align}
Using \eqr{e:abd} and $|A|^2 \leq 4$, we get
\begin{equation} \label{e:stab3}
\int |A|^2 \psi^2
\leq 4 \, C_1 \, R^3 + \int |A|^2 \chi^2 \psi^2
\leq 10 \, C_1 \, R^3 + 3 \, R^{-2} \, \Area \, (\cB_R) \, .
\end{equation}
The lemma follows
from \eqr{e:stab3} and Corollary \ref{c:lcacc}.
\end{proof}

The polynomial growth allows us to find large intrinsic balls
with a fixed doubling:

\begin{Cor} \label{c:choosing}
There exists $C_2$ so that
given $\beta , R_0 > 1$, we get $R_{2}$
so: If $0 \in \Sigma \subset B_{R_{2}}$ is an embedded
minimal disk, $\partial \Sigma \subset \partial B_{R_2}$,
$|A|^2(0) = 1$, and $|A|^2 \leq 4$,
then there exists $R_0 < R < R_2 / (2 \, \beta)$ with
\begin{equation}
\int_{\cB_{3 \, R}} |A|^2 +
\beta^{-10} \, \int_{\cB_{2 \, \beta \, R}} |A|^2 \leq C_{2} \, R^{-2} \,
\Area \, ( \cB_{R} ) \, . \label{e:fe3}
\end{equation}
\end{Cor}

\begin{proof}
Set $\ca (s) = \Area (\cB_s)$.
Given $m$,
Lemma \ref{l:vgr} gives
\begin{equation} \label{e:fe1}
\left( \min_{1 \leq n \leq m} \frac{\ca ( (4\beta)^{2n} \, R_0)}
{\ca ( (4\beta)^{2n-2} \, R_0)} \right)^m \leq
\frac{\ca ((4\beta)^{2m} \, R_0 )}{\ca (R_0)} \leq C_1' \,
(4\beta)^{10m} \, R_0^3 \, .
\end{equation}
Fix $m$ with $C_1' \, R_0^3 < 2^{m}$ and
set
$R_2 = 2 \, (4\beta)^{2m} \, R_0$. By \eqr{e:fe1}, there exists
$R_1 = (4\beta)^{2n-2} \, R_0$ with $1 \leq n \leq m$ so
\begin{equation} \label{e:fe2}
\frac{\ca ( (4\beta)^{2} \, R_1)}
{\ca ( R_1)}
\leq 2 \, (4 \beta)^{10} \, .
\end{equation}
For simplicity, assume that $\beta = 4^q$ for $q \in \ZZ^+$.
As in \eqr{e:fe1}, \eqr{e:fe2}, we get $0 \leq j \leq q$ with
\begin{equation} \label{e:fe2a}
\frac{\ca ( 4^{j+1} \, R_1)}
{\ca ( 4^{j} \, R_1)}
\leq
\left[ \frac{\ca ( 4 \beta \, R_1)}
{\ca ( R_1)} \right]^{1/(q+1)}
\leq
2^{1/(q+1)} \, 4^{10} \, .
\end{equation}
Set $R = 4^{j} \, R_1$.
Combining \eqr{e:fe2}, \eqr{e:fe2a}, and
Corollary \ref{c:lcacc} gives
\eqr{e:fe3}.
\end{proof}

\section{The local structure near the axis} \label{s:s2}

\begin{figure}[htbp]
    \setlength{\captionindent}{20pt}
    \begin{minipage}[t]{0.5\textwidth}
\centering\input{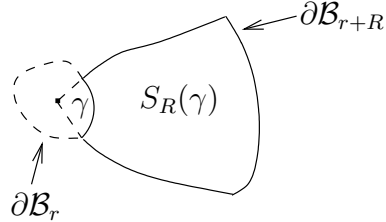}
    \caption{The intrinsic sector over a curve $\gamma$
    defined in \eqr{e:defsip}.}
\label{f:9}
\end{minipage}
\end{figure}

 Given $\gamma \subset  \partial \cB_{r}$,  define the intrinsic
sector, see fig. \ref{f:9},
\begin{equation} \label{e:defsip}
  S_{R}(\gamma) = \{ \exp_0 (v) \, | \,  r \leq |v| \leq  r + R {\text{ and }}
\exp_0 (r \, v / |v|) \in   \gamma \} \, .
\end{equation}

The key for proving  Theorem \ref{t:blowupwindinga} is to
find $n$ large intrinsic sectors with a scale-invariant
curvature bound.
To do this,
 we first
use Corollary \ref{c:scsi}  to  bound
${\text{Length}}(\partial  \cB_R)/R$ from below for $R \geq R_0$.
Corollary \ref{c:choosing} gives  $R_3 >R_0$
and $n$ long disjoint curves $\tilde{\gamma}_i \subset \partial \cB_{R_3}$
so the sectors over $\tilde{\gamma}_i$ have bounded
$\int |A|^2$.  Corollary \ref{c:scsi}  gives the
 curvature  bound.
Once we have these sectors, for $n$ large,
two  must be
close and hence, by Lemmas  \ref{l:fiscm} and \ref{l:dnl},
$1/2$-stable.  The $N$-valued graph is then given
by corollary II.1.34 of \cite{CM3}:

\begin{Cor} \label{c:uselater}
\cite{CM3}.
Given $\omega > 8, 1 > \epsilon > 0, C_0$, and $N$, there exist $m_1 ,
\Omega_1$ so: If $0 \in \Sigma$ is an embedded minimal disk,
$\gamma \subset \partial \cB_{r_1}$ is a
curve, $\int_{\gamma} k_g  < C_0 \, m_1$,   ${\text{Length}}(\gamma) =
m_1  \, r_1$, and $\cT_{r_1 / 8} (S_{\Omega_1^2 \, \omega \, r_1 } (\gamma))$
is $1/2$-stable,    then
(after rotating $\RR^3$) $S_{ \Omega_1^2 \, \omega \, r_1}
(\gamma)$ contains an $N$-valued graph $\Sigma_N$ over $D_{\omega \,
\Omega_1 \, r_1} \setminus D_{\Omega_1 \, r_1}$ with gradient
$\leq \epsilon$, $|A| \leq \epsilon /  r$, and
$\dist_{S_{ \Omega_1^2 \, \omega \, r_1}
(\gamma)} ( \gamma , \Sigma_N ) < 4 \, \Omega_1 \, r_1$.
\end{Cor}

\begin{figure}[htbp]
    \setlength{\captionindent}{20pt}
    \begin{minipage}[t]{0.5\textwidth}
\centering\input{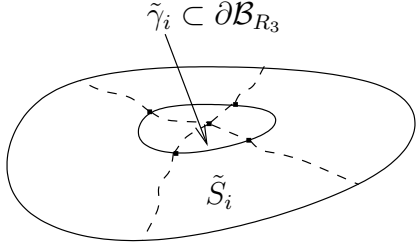}
    \caption{Equation \eqr{e:defsiptilde} divides a punctured
    ball into
    sectors $\tilde S_i$.}
\label{f:10}
\end{minipage}
\end{figure}

\begin{proof}
(of Theorem \ref{t:blowupwindinga}).
Rescale by $C/r_0$ so that $|A|^2(0)=1$ and
$|A|^2 \leq 4$ on $B_C$.

Let $C_2$ be from Corollary \ref{c:choosing} and then let
$m_1 , \Omega_1 > \pi$ be
given by Corollary \ref{c:uselater} with $C_0$ there
$=2 \, C_2 + 2$.  Fix $a_0$ large (to be
chosen). By Corollaries \ref{c:lcacc}, \ref{c:scsi},
there exists $R_0 = R_0(a_0)$ so
that for any $R_3 \geq R_0$
\begin{equation} \label{e:a5}
  a_0 \, R_3 \leq R_3/4 \, \int_{\cB_{R_3/2}} |A|^2  \leq
{\text{Length}} (\partial \cB_{R_3})
   \, .
\end{equation}
Set $\beta = 2 \Omega_1^2 \, \omega$.
Corollaries \ref{c:lcacc}, \ref{c:choosing}  give $R_2 =
R_2 (R_0 , \beta)$ so if $C \geq R_2$, then
there is $R_0 < {R_3} < R_2 / (2 \, \beta)$ with
\begin{equation}
\int_{\cB_{3 \, {R_3}}} |A|^2 +
\beta^{-10} \, \int_{\cB_{2 \, \beta \, {R_3}}} |A|^2 \leq
C_2  \, {R_3}^{-2} \, \Area (\cB_{R_3} ) \leq
C_2 \, {\text{Length}} (\partial \cB_{R_3})
/ (2{R_3})
 \, . \label{e:a4}
\end{equation}
Using \eqr{e:a5}, choose $n$ so that
\begin{equation} \label{e:fixn}
a_0 \, {R_3} \leq 4 \, m_1 \, n \, {R_3} < {\text{Length}} (\partial \cB_{R_3})
 \leq 8 \, m_1 \, n \, {R_3} \, ,
\end{equation}
and fix $2n$ disjoint curves $\tilde \gamma_i \subset
\partial \cB_{{R_3}}$ with length $2 \, m_1 \, {R_3}$.
Define the intrinsic  sectors (see fig. \ref{f:10})
\begin{equation} \label{e:defsiptilde}
\tilde S_i = \{ \exp_0 (v) \, | \, 0<|v| \leq 2 \, \beta \, {R_3} {\text{ and }}
\exp_0 ({R_3} \, v / |v|) \in \tilde  \gamma_i \} \, .
\end{equation}
Since the $\tilde S_i$'s are disjoint, combining \eqr{e:a4} and \eqr{e:fixn} gives
\begin{equation} \label{e:a6}
\sum_{i=1}^{2n} \left( \int_{\cB_{3\,{R_3}} \cap \tilde S_i} |A|^2 + \beta^{-10} \,
\int_{\tilde S_i} |A|^2 \right) \leq 4 \, C_2 \, m_1 \, n \, .
\end{equation}
Hence, after reordering the $\tilde \gamma_i$, we can assume that for
$1 \leq i \leq n$
\begin{equation} \label{e:a7a}
\int_{\cB_{3\,{R_3}} \cap \tilde S_i} |A|^2 + \beta^{-10} \,
\int_{\tilde S_i} |A|^2 \leq 4 \, C_2 \, m_1 \, .
\end{equation}
Using the Riccati comparison theorem, there are curves $\gamma_i \subset
\partial \cB_{2{R_3}} \cap \tilde  S_i$ with length $2 \, m_1 \, {R_3}$ so that if
$y \in S_i=S_{\beta R_3}(\gamma_i) \subset \tilde S_i$,
then
$\cB_{\dist_{\Sigma}(0,y)/2}(y) \subset \tilde S_i$.
Hence, by Corollary
\ref{c:scsi} and
\eqr{e:a7a}, we get for $y \in S_i$ and $i\leq n$
\begin{equation} \label{e:a9}
\sup_{\cB_{\dist_{\Sigma}(0,y)/4}(y)} |A|^2 \leq
C_3 \, \dist_{\Sigma}^{-2} (0,y)  \, ,
\end{equation}
where $C_3 = C_3 (\beta ,  m_1 )$. For $i \leq n$,
\eqr{e:a7a} and the Gauss-Bonnet theorem yield
\begin{equation} \label{e:a7b}
\int_{\gamma_i} k_g \leq 2 \, \pi + 2 \, C_2 \, m_1 < (2\, C_2 + 2) \, m_1 \, .
\end{equation}
By \eqr{e:a9} and a Riccati comparison argument,
there exists $C_4 = C_4 (C_3)$ so that for $i\leq n$
\begin{equation} \label{e:a7bb}
1 / (2{R_3}) \leq \min_{\gamma_i} k_g \leq \max_{\gamma_i} k_g \leq C_4 / {R_3} \, .
\end{equation}

Applying Lemma \ref{l:dnl} repeatedly (and using \eqr{e:a9}),
it is easy to see that there exists $\alpha > 0$
so that if $i_1 < i_2 \leq n$ and
\begin{equation} \label{e:getcloseA}
\dist_{C^{1}([0,2m_1],\RR^3)}
(\gamma_{i_1}/{R_3} , \gamma_{i_2}/{R_3} )
\leq \alpha \, ,
\end{equation}
then
$\{ z + u(z) \, \nn(z) \, | \, z \in \cT_{{R_3}/4} (S_{i_1}) \} \subset
\cup_{y \in S_{i_2}} \cB_{\dist_{\Sigma}(0,y)/4}(y)$
for a function $u \ne 0$ with
\begin{equation} \label{e:graphca}
|\nabla u| + |A| \, |u| \leq C_0' \, \dist_{C^{1}([0,2m_1],\RR^3)}
(\gamma_{i_1}/{R_3} , \gamma_{i_2}/{R_3} ) \, \, .
\end{equation}
Here $\dist_{C^{1}([0,2m_1],\RR^3)} (\gamma_{i_1}/{R_3} , \gamma_{i_2}/{R_3} )$
is the scale-invariant $C^1$-distance between the
curves.

Next, we use compactness to show that
\eqr{e:getcloseA} must hold for $n$ large.
Namely, since each
$\gamma_i / {R_3} \subset B_{2}$ is parametrized by arclength on
$[0,2 m_1]$ and has a uniform
$C^{1,1}$ bound
by \eqr{e:a7bb}, this set of maps is compact by
the Arzela-Ascoli theorem. Hence, there exists $n_0$
 so that
if $n \geq n_0$, then \eqr{e:getcloseA} holds
for some $i_1 < i_2 \leq n$. In particular,
\eqr{e:graphca} and
Lemma \ref{l:fiscm} imply that $S_{i_1}$ is
$1/2$-stable for $n$ large (now choose $a_0, R_0 , R_2$).
After rotating $\RR^3$, Corollary \ref{c:uselater} gives
the $N$-valued graph
$\Sigma_g \subset S_{i_1}$ over $D_{2 \omega \, \Omega_1 \, {R_3}} \setminus
D_{ 2 \Omega_1 \, {R_3}}$ with gradient $\leq \epsilon$, $|A| \leq \epsilon / r$,
and $\dist_{\Sigma} (0 , \Sigma_g) \leq 8 \, \Omega_1 \, {R_3}$.
Rescaling by $r_0/C$,
the theorem follows with $\bar{R} = 2 \Omega_1 \, {R_3} r_0 / C$.
\end{proof}

\begin{Cor} \label{c:blowupwinding}
Given $N > 1$ and $\tau > 0$, there exist $\Omega > 1$ and
$C > 0$ so: Let
$0\in \Sigma^2\subset B_{R}$ be an embedded minimal
disk, $\partial \Sigma\subset \partial B_{R}$. If
$R>r_0>0$ with
$\sup_{B_{r_0} \cap \Sigma}|A|^2\leq 4\,C^2\,r_0^{-2}$
and $|A|^2(0)=C^2\,r_0^{-2}$, then there exists
(after a rotation)
an $N$-valued graph $\Sigma_g \subset \Sigma$ over $D_{R/\Omega}
\setminus D_{r_0}$ with gradient $\leq \tau$,
$\dist_{\Sigma}(0,\Sigma_g) \leq 4 \, r_0$, and
$\Sigma_g \subset \{ x_3^2 \leq \tau^2 \, (x_1^2 + x_2^2) \}$.
\end{Cor}

\begin{proof}
This follows immediately by combining
Theorems \ref{t:blowupwindinga} and \ref{t:spin4ever2}.
\end{proof}

\begin{figure}[htbp]
    \setlength{\captionindent}{20pt}
    \begin{minipage}[t]{0.5\textwidth}
\centering\input{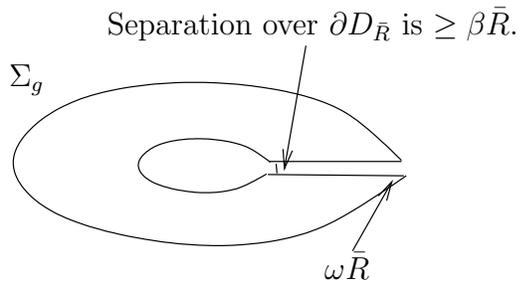}
    \caption{Proposition \ref{p:blowupgap}:
    The initial separation is inversely proportional
    to the maximum of $|A|$.}
\label{f:11}
\end{minipage}
\end{figure}

\begin{Pro} \label{p:blowupgap}
See fig. \ref{f:11}.
There exists $\beta > 0$ so:
If $\Sigma_g \subset \Sigma$ is as in
Theorem \ref{t:blowupwindinga},
then the separation between the sheets of
$\Sigma_g$ over $\partial D_{\bar{R}}$ is at least $\beta \, \bar{R}$.
\end{Pro}

\begin{proof}
This follows easily from the curvature bound,
Lemma \ref{l:dnl},
the Harnack inequality, and
estimates for $1/2$-stable surfaces.
\end{proof}

\section{The blow up}

Combining Corollary \ref{c:blowupwinding} and a blowup argument
will give Theorem \ref{t:blowupwinding0}.

\begin{Lem} \label{l:bup}
If $0 \in \Sigma \subset B_{r_0}$, $\partial \Sigma \subset \partial
B_{r_0}$, and
$\sup_{B_{r_0/2} \cap \Sigma}|A|^2\geq 16\,C^2\,r_0^{-2}$, then there exist
$y \in \Sigma$ and $r_1 < r_0 - |y|$ with $|A|^2 (y)
= C^2 \,r_1^{-2}$ and
$\sup_{B_{r_1}(y) \cap \Sigma}|A|^2\leq 4 \, C^2 \,r_1^{-2}$.
\end{Lem}

\begin{proof}
Set $F(x)=(r_0-|x|)^2\, |A|^2(x)$.
Since $F \geq 0$, $F|\partial B_{r_0}\cap \Sigma=0$, and
$\Sigma$ is compact,
$F$ achieves its maximum at $y \in \partial
B_{r_0-\sigma}\cap \Sigma$ with $0< \sigma \leq r_0$.
Since
$\sup_{B_{r_0/2} \cap \Sigma}|A|^2\geq 16 \,C^2\,r_0^{-2}$,
\begin{equation} \label{e:o2.0a}
F(y)=\sup_{B_{r_0}\cap \Sigma} F \geq
4 \, C^2 \, .
\end{equation}
To get the first claim, define $r_1>0$ by
\begin{equation} \label{e:o2.1a}
r_1^2\,|A(y)|^2= C^2\, .
\end{equation}
Since $F(y)=\sigma^2 \, |A(y)|^2 \geq 4 \, C^2$, we have
$2 \, r_1 \leq \sigma$. Finally, by \eqr{e:o2.0a},
\begin{equation} \label{e:o2.3a}
\sup_{B_{r_1}(y)\cap \Sigma}
\left( \frac{\sigma}{2} \right) ^2\, |A|^2
\leq \sup_{B_{\frac{\sigma}{2}}(y)\cap \Sigma}
\left( \frac{\sigma}{2} \right) ^2\, |A|^2
\leq \sup_{B_{\frac{\sigma}{2}}(y)\cap \Sigma} F
\leq \sigma^2\, |A(y)|^2\, .
\end{equation}
\end{proof}

\begin{proof}
(of Theorem \ref{t:blowupwinding0}).
This follows immediately from Corollary \ref{c:blowupwinding} and
Lemma \ref{l:bup}.
\end{proof}

If
$y_i\in \Sigma_i$ is a sequence of minimal disks with $y_i \to y$ and $|A|(y_i)$
blowing up, then we can take $r_0 \to 0$ in
Theorem \ref{t:blowupwinding0}. Combining this with the sublinear
growth of the separation between the sheets from
\cite{CM3}, we will get in Theorem \ref{t:stablim}
a smooth limit through $y$.

Below $\Sigma^{0,2\pi}_{r,s} \subset \Sigma$ is the ``middle sheet'' over
$\{ (\rho, \theta) \, | \, 0 \leq \theta
\leq 2 \pi , \, r \leq \rho \leq s \}$.
The sublinear growth is given by proposition II.2.12 of \cite{CM3}:

\begin{Pro} \label{l:grades1}
\cite{CM3}.  See fig. \ref{f:12}.
Given $\alpha > 0$, there exist $\delta_p > 0 , N_g > 5 $ so:
If $\Sigma$ is a $N_g$-valued minimal graph over
$D_{\e^{N_g} \, R} \setminus D_{\e^{-N_g} \, R}$ with gradient
$\leq 1$ and $0 < u < \delta_p \, R$ is a solution of the minimal graph
equation over $\Sigma$ with $|\nabla u| \leq 1$,
then for $R \leq s \leq 2 \, R$
\begin{align} \label{e:wantit}
\sup_{ \sztp_{R,2R} } |A_{\Sigma}| &+
\sup_{ \sztp_{R,2R} } |\nabla u| / u \leq
\alpha / (4\,R) \, , \\ \label{e:slg}
\sup_{ \sztp_{R,s} } u &\leq (s/R)^{\alpha} \, \sup_{
\sztp_{R,R} } u \, .
\end{align}
\end{Pro}

\begin{Thm} \label{t:stablim}
See fig. \ref{f:13}.
There exists $\Omega > 1$ so:
Let $y_i\in \Sigma_i \subset B_{R}$ with $\partial \Sigma_i \subset
\partial B_{R}$ be embedded minimal disks where $y_i \to 0$. If $|A_{\Sigma_i}|(y_i)
\to \infty$, then, after a rotation and passing to a subsequence,
there exist
$\epsilon_i \to 0$, $\delta_i \to 0$, and $2$-valued minimal
graphs $\Sigma_{d,i} \subset \{ x_3^2 \leq x_1^2 + x_2^2 \} \cap
\Sigma_i$ over $D_{R/\Omega} \setminus D_{\epsilon_i}$
with gradient $\leq 1$, and separation at most $\delta_i \, s$ over
$\partial D_s$. Finally, the $\Sigma_{d,i}$ converge (with
multiplicity two) to a smooth minimal graph through $0$.
\end{Thm}

\begin{proof}
The first part follows immediately from
taking $r_0 \to 0$ in
Theorem \ref{t:blowupwinding0}. When $s$ is small, the bound
on the separation follows from the gradient bound.
The separation then grows less than linearly
by Proposition \ref{l:grades1}, giving the bound for large
$s$ and showing that the $\Sigma_{d,i}$ close up in the limit.
In particular, the $\Sigma_{d,i}$ converge to a minimal graph
$\Sigma'$ over $D_{R/\Omega} \setminus \{ 0 \}$ with gradient $\leq 1$
and $\Sigma' \subset \{ x_3^2 \leq x_1^2 + x_2^2 \}$. By
a standard removable singularity theorem,
$\Sigma' \cup \{ 0 \}$ is a smooth minimal graph over
$D_{R/\Omega}$.
\end{proof}

\begin{figure}[htbp]
    \setlength{\captionindent}{20pt}
    \begin{minipage}[t]{0.5\textwidth}
    \centering\input{shn7.pstex_t}
    \caption{The sublinear growth of the separation $u$ of
the multi-valued graph $\Sigma$:
    $u(2R) \leq 2^{\alpha} \, u(R) $ with $\alpha < 1$.}
    \label{f:12}
\end{minipage}\begin{minipage}[t]{0.5\textwidth}
    \centerline{\input{blow13.pstex_t}}
 \caption{Theorem \ref{t:stablim}: As $|A_{\Sigma_i}|(y_i) \to \infty$ and
    $y_i \to y$, $2$-valued
    graphs converge to a graph through $y$.
(The upper sheets of the $2$-valued graphs collapses to the lower sheets.)}
\label{f:13}
\end{minipage}
\end{figure}

\appendix

\end{document}